\documentclass[twoside             
]{amsart}
\usepackage[centertags]{amsmath}

\usepackage{amsfonts}
\usepackage{amssymb}
\usepackage{amsthm}
\usepackage{eucal}
\usepackage[cmtip,all]{xy}
\usepackage{graphicx}

\newtheorem{thm}[equation]{Theorem}

\newtheorem{lem}[equation]{Lemma}

\theoremstyle{definition}
\newtheorem{defn}[equation]{Definition}
\theoremstyle{remark}
\newtheorem{rem}[equation]{Remark}
\newtheorem{exm}[equation]{Example}
\numberwithin{equation}{section}

\newcommand{\set}[1]{\left\{#1\right\}}

\newcommand{\To}{\longrightarrow}

\newcommand{\A}{\mathcal{A}}

\newcommand{\C}[1]{\mathcal{#1}}
\newcommand{\CC}[1]{\mathbf{#1}}
\newcommand{\D}[1]{\mathcal{#1}}
\newcommand{\DD}[1]{\mathcal{#1}^{\mathrm{der}}}
\newcommand{\der}{\mathbb{D}}
\newcommand{\db}{D^b}
\newcommand{\cb}{C^b}
\newcommand{\hb}{H^b}

\newcommand{\relnn}{{\textrm{\rm R}}}
\newcommand{\relss}{{\textrm{\rm DR}}}
\newcommand{\reln}[1]{{\textrm{\rm(\relnn #1)}}}
\newcommand{\rels}[1]{{\textrm{\rm(\relss #1)}}}

\newcommand{\gennn}{{\textrm{\rm G}}}
\newcommand{\genss}{{\textrm{\rm DG}}}
\newcommand{\genn}[1]{{\textrm{\rm(\gennn #1)}}}
\newcommand{\gens}[1]{{\textrm{\rm(\genss #1)}}}

\newcommand{\EZDIAG}[5]{\xymatrix@C+=1.5cm{*+[r]{#1}
\ar@(u,l)_(0.62){\displaystyle #5}[]
\ar@<.7ex>^-{#3}[r]&\ar@<.7ex>^-{#4}[l]#2}}

\def\r{\rightarrow} 

\def\into{\rightarrowtail}
\def\onto{\twoheadrightarrow}

\def\estrella{}

\def\ner{\operatorname{Ner}}

\def\ho{\operatorname{Ho}}

\def\diag{\operatorname{Diag}}

\def\st{\stackrel} 

\renewcommand{\ker}{\operatorname{Ker}}

\def\Z{\mathbb{Z}}

\newcommand{\grupo}[1]{\langle #1\rangle}

\newcommand{\bak}[1]{#1'}

\begin{document}

\title{Maltsiniotis's first conjecture for $K_1$}%
\author{Fernando Muro}%
\address{Universitat de Barcelona, Departament d'\`Algebra i Geometria, Gran via de les corts catalanes 585, 08007 Barcelona, Spain}
\email{fmuro@ub.edu}

\thanks{The author was partially supported
by the Spanish Ministry of Education and Science under MEC-FEDER grants MTM2004-03629 and MTM2007-63277, and a Juan de la Cierva research contract.}

\begin{abstract}
We show that $K_1(\C{E})$ of an exact category $\C{E}$ agrees with $K_1(\der{\C{E}})$ of the associated triangulated derivator $\der{\C{E}}$. More generally we show that $K_1(\C{W})$ of a Waldhausen category $\C{W}$ with cylinders and a saturated class of weak equivalences agrees with $K_1(\der{\C{W}})$ of the associated right pointed derivator $\der{\C{W}}$.
\end{abstract} \maketitle 

\section*{Introduction}

\renewcommand{\theequation}{\Alph{equation}}

For a long time there was an interest in defining a nice $K$-theory for triangulated categories such that Quillen's $K$-theory of an exact category $\C{E}$ agrees with the $K$-theory of its bounded derived category $\db(\C{E})$. Schlichting \cite{kttc} showed that such a $K$-theory for triangulated categories cannot exist. It was then natural to ask about the definition of a nice $K$-theory for algebraic structures interpolating between $\C{E}$ and $\db(\C{E})$. 

The best known intermediate structure is $\cb(\C{E})$, the Waldhausen category of bounded complexes in $\C{E}$, with quasi-isomorphisms as weak equivalences and cofibrations given by chain morphisms which are levelwise admissible monomorphisms. The derived category $\db(\C{E})$ is the localization of $\cb(\C{E})$ with respect to weak equivalences. The Gillet-Waldhausen theorem\footnote{The proof due to Thomason-Trobaugh \cite[Theorem 1.11.7]{tt} corrects Gillet's \cite[6.2]{rrthakt} and uses an extra hypothesis on $\C{E}$. This hypothesis is not strictly necessary, since the general case follows then from cofinality arguments, see \cite{tckt}.}, relating Quillen's $K$-theory to Waldhausen's $K$-theory, states that the homomorphisms
$$\tau_n\colon K_n(\C{E})\To K_n(\cb(\C{E})),\quad n\geq0,$$ induced by 
the inclusion $\C{E}\subset\cb(\C{E})$ of complexes concentrated in degree $0$, are isomorphisms. 

The category $\cb(\C{E})$ is considered to be too close to $\C{E}$ so one would still like to find an algebraic stucture with a nice $K$-theory interpolating between $\cb(\C{E})$ and $\db(\C{E})$. The notion of a triangulated derivator \cite{derivateurs, ktdt} seems to be a strong candidate. 

Maltsiniotis \cite{ktdt} defined a $K$-theory for triangulated derivators together with natural homomorphisms
$$\rho_n\colon K_n(\C{E})\To K_n(\der{\C{E}}),\quad n\geq0,$$
where $\der{\C{E}}$ is the triangulated derivator associated to an exact category $\C{E}$, constructed by Keller in the appendix of \cite{ktdt}. Cisinski and Neeman proved the additivity of triangulated derivator $K$-theory \cite{adkt}. Maltsiniotis also conjectured that $\rho_n$ is an isomorphism for all $n$. He succeeded in proving the conjecture for $n=0$. 

The following theorem is the main result of this paper.

\begin{thm}\label{main}
Let $\C{E}$ be an exact category. The natural homomorphism $$\rho_1\colon K_1(\C{E})\st{\cong}\To K_1(\der{\C{E}})$$
is an isomorphism.
\end{thm}



In order to obtain Theorem \ref{main} we use techniques introduced in \cite{1tk}. There we give a presentation of an abelian $2$-group $\D{D}_*\C{W}$ which encodes $K_0(\C{W})$ and $K_1(\C{W})$ of a Waldhausen category $\C{W}$, and moreover the $1$-type of the $K$-theory spectrum $K(\C{W})$ whose homotopy groups are the $K$-theory groups of $\C{W}$. This presentation is a higher dimensional analogue of the classical presentation of $K_0(\C{W})$. Here we similarly define an abelian $2$-group $\DD{D}_*\C{W}$ which models the $1$-type of the $K$-theory spectrum $K(\der\C{W})$ of the right\footnote{The references \cite{derivateurs,ciscd} and \cite{sdckt1,cht} follow a different convention with respect to sides. Here we follow the convention in \cite{derivateurs, ciscd}, so what we call a `right pointed derivator' is the same as a `left pointed derivator' in \cite{sdckt1}.} pointed derivator $\der\C{W}$ associated to a Waldhausen category $\C{W}$ with cylinders and a saturated class of weak equivalences, such as $\C{W}=\cb(\C{E})$. The $K$-theory for this kind of derivators, more general than triangulated derivators, was  defined by Garkusha \cite{sdckt1} extending the work of Maltsiniotis \cite{ktdt}. There are defined comparison homomorphisms
$$\mu_n\colon K_n(\C{W})\To K_n(\der\C{W}),\quad n\geq0.$$
These homomorphisms cannot be isomorphisms in general, as shown in \cite{ktsc}. Nevertheless
we here prove the following result.

\begin{thm}\label{main2}
Let $\C{W}$ be a Waldhausen category with cylinders and a saturated class of weak equivalences. The natural homomorphism $$\begin{array}{c}\mu_0\colon K_0(\C{W})\st{\cong}\To K_0(\der\C{W}),\\
\mu_1\colon K_1(\C{W})\st{\cong}\To K_1(\der\C{W}),
\end{array}
$$ 
are isomorphisms.
\end{thm}

In Remark \ref{last} we comment on the case where the hypothesis on the saturation of weak equivalences is replaced by the 2 out of 3 axiom, which is a weaker assumption.

Theorem \ref{main} is actually a corollary of the Gillet-Waldhausen theorem and Theorem \ref{main2}, since 
$\der\cb(\C{E})=\der\C{E}$ and the natural homomorphisms $\rho_n$ factor as
$$\rho_n\colon K_n(\C{E})\st{\tau_n}\To K_n(\cb(\C{E}))\st{\mu_n}\To K_n(\der{\C{E}}),\quad n\geq0.$$

We assume the reader certain familiarity with exact, Waldhausen and derived categories, with simplicial constructions and with homotopy theory. We refer to \cite{wiak, mha, gj} for the basics.

\subsection*{Acknowledgements}

I am very grateful to Grigory Garkusha for suggesting the possibility of using \cite{1tk} in order to tackle Maltsiniotis's first conjecture in dimension $1$. I also feel indebted to Denis-Charles Cisinski, who kindly indicated how to extend the results of a preliminary version of this paper to a broader generality.

\renewcommand{\theequation}{\thesection.\arabic{equation}}

\section{The bounded derived category of an exact category}\label{ho}

In this section we outline the two-step construction of the derived category $\db(\C{E})$ of an exact category $\C{E}$. This construction is a special case of the homotopy category $\ho \C{W}$ of a Waldhausen category $\C{W}$ with cylinders satisfying the 2 out of 3 axiom, $\db(\C{E})=\ho\cb(\C{E})$.

\begin{defn}
A \emph{Waldhausen category} is a category $\C{W}$ with a distinguished zero object $0$ and two distinguished subcategories $w\C{W}$ and $c\C{W}$, whose morphisms are called \emph{cofibrations} and \emph{weak equivalences}, respectively.  A morphism which is both a weak equivalence and a cofibration is said to be a \emph{trivial cofibration}.
The arrow $\rightarrowtail$ stands for a cofibration and $\st{\sim}\r$ for a weak equivalence. 
\begin{itemize}
\item All morphisms $0\r A$ are cofibrations. All isomorphisms are cofibrations and weak equivalences. 

\item The push-out of a morphism along a cofibration is always defined
$$\xymatrix{A\ar@{>->}[r]\ar[d]\ar@{}[rd]|{\text{push}}&B\ar[d]\\X\ar@{>->}[r]&X\cup_AB}$$
and the lower map is also a cofibration. 

\item Given a commutative diagram
$$\xymatrix{X\ar[d]_\sim&A\ar[d]_\sim\ar[l]\ar@{>->}[r]&B\ar[d]_\sim\\X'&A'\ar[l]\ar@{>->}[r]&B'}$$
the induced map $X\cup_AB\st{\sim}\r X'\cup_{A'}B'$ is a weak equivalence. 
\end{itemize}
Notice that coproducts $A\vee B=A\cup_0B$ are defined in $\C{W}$. 

A functor $\C{W}\r\C{W}'$ between Waldhausen categories is \emph{exact} if it preserves cofibrations, weak equivalences, push-outs along cofibrations and the distinguished zero object.
\end{defn}

\begin{exm}
Recall that an \emph{exact category} $\C{E}$ is a full subcategory of an abelian category $\C{A}$ such that $\C{E}$ contains a zero object of $\C{A}$ and $\C{E}$ is closed under extensions in $\C{A}$. A \emph{short exact sequence} in $\C{E}$ is a short exact sequence in $\C{A}$ between objects in $\C{E}$. A morphism in $\C{E}$ is an \emph{admissible monomorphism} if it is the initial morphism of some short exact sequence.
The category $\C{E}$ is a Waldhausen category with admissible monomorphisms as cofibrations and isomorphisms as weak equivalences. In order to complete the structure we fix a zero object $0$ in $\C{E}$.

We denote by $\cb(\C{E})$ the category of bounded complexes in $\C{E}$,
$$\cdots \r A^{n-1}\st{d}\To A^n\st{d}\To A^{n+1}\r\cdots, \quad d^2=0, \quad A^n=0\text{ for }|n|\gg0.$$
A chain morphism $f\colon A^{\estrella}\r B^{\estrella}$ in $\cb(\C{E})$ is a \emph{quasi-isomorphism} if it induces an isomorphism in homology computed in the ambient abelian category $\C{A}$. The category $\cb(\C{E})$ is a Waldhausen category. Weak equivalences are quasi-isomorphisms and cofibrations are levelwise admissible monomorphisms. 

There is a full exact inclusion of Waldhausen categories $\C{E}\subset\cb(\C{E})$ sending an object $X$ in $\C{E}$ to the complex
$$\cdots\r 0\r X\r 0\r\cdots,$$
with $X$ in degree $0$. 
\end{exm}

\begin{defn}\label{hotsat}
The \emph{homotopy category} $\ho\C{W}$ of a Waldhausen category is a category equipped with a functor
$$\zeta\colon \C{W}\To \ho\C{W}$$
sending weak equivalences to isomorphisms. Moreover, $\zeta$ is initial among all functors $\C{W}\r\C{C}$ sending weak equivalences to isomorphisms, so $\ho\C{W}$ is well defned up to canonical isomorphism over $\C{W}$. This category can be constructued as a category of fractions, in the sense of \cite{gz}, by formally inverting weak equivalences in $\C{W}$. 

The class of weak equivalences is \emph{saturated} if any morphism $f\colon A\r B$ in $\C{W}$ such that $\zeta(f)$ is an isomorphism in $\ho\C{W}$ is indeed a weak equivalence $f\colon A\st{\sim}\r B$.
\end{defn}

\begin{exm}
Weak equivalences in $\cb(\C{E})$, i.e. quasi-isomorphisms, are saturated since they are detected by a functor $H^*\colon\cb(\C{E})\r\C{A}^\Z$, the cohomology functor from bounded complexes in $\C{E}$ to $\Z$-graded objects in $\C{A}$, see \cite[Proposition 1.1]{scm}.
\end{exm}

The homotopy category always exists up to set theoretical difficulties which do not arise if $\C{W}$ is a small category, for instance. This is not a harmful assumption if one is interested in $K$-theory since smallness may also be required in order to have well defined $K$-theory groups. The homotopy category can however be constructed in a more straighforward way if the Waldhausen category $\C{W}$ satisfies further properties.

\begin{defn}
A Waldhausen category $\C{W}$ satisifies the \emph{2 out of 3 axiom} provided given a commutative diagram in $\C{W}$
$$\xymatrix@C=20pt{&C&\\A\ar[ru]\ar[rr]&&B\ar[lu]}$$
if two arrows are weak equivalences then the third one is also a weak equivalence.


Given an object $A$ in $\C{W}$ a \emph{cylinder} $IA$ is an object together with a factorization
of the folding map $(1,1)\colon A\vee A\r A$ as a cofibration followed by a weak equivalence,
$$A\vee A\mathop{\rightarrowtail}\limits_i IA\mathop{\To}\limits^\sim_p A.$$ We say that $\C{W}$ has cylinders if all objects have a cylinder.
\end{defn}

\begin{exm}
The Waldhausen category $\cb(\C{E})$ has cylinders. The cylinder of a bounded complex $A$ can be functorially chosen as
$$(IA)^n\;=\;A^n\oplus A^{n+1}\oplus A^n,\quad d=\left(\begin{array}{ccc}d&-1&0\\0&-d&0\\0&1&d\end{array}\right)\colon (IA)^n\To (IA)^{n+1}.$$
\end{exm}

\begin{rem}
The 2 out of 3 axiom is often called the saturation axiom. We do not use this terminology in this paper in order to avoid confusion with Definition \ref{hotsat}. 

Usually one considers more structured cylinders in Waldhausen categories, compare \cite[Definition IV.6.8]{wiak}. For the purposes of this paper it is enough to consider cylinders as defined above.
\end{rem}

\begin{rem}
As one can easily check, a Waldhausen category with a saturated class of weak equivalences satisfies the 2 out of 3 axiom. This applies to $\cb(\C{E})$.

A Waldhausen category with cylinders $\C{W}$ satisfying the 2 out of 3 axiom is an example of a \emph{right derivable category}, in the sense of \cite{ciscd}, also called \emph{precofibration category} in \cite{cht}, see \cite[Example 2.23]{ciscd} or \cite[Proposition 2.4.2]{cht}. In particular any morphism in $\C{W}$ can be factored as a cofibration followed by a weak equivalence which is left inverse to a trivial cofibration, see \cite[Proposition 1.3.1]{cht}. Moreover, one can define a homotopy relation in $\C{W}$ and construct the homotopy category $\ho\C{W}$ by a homotopy calculus of left fractions as we indicate below, see \cite[Section 1]{ciscd} or \cite[Section 5.4]{cht}.
\end{rem}

Let $\C{W}$ be a Waldhausen category with cylinders satisfying the 2 out of 3 axiom. As usual we say that two morphisms $f,g\colon A\r B$ in $\C{W}$ are \emph{strictly homotopic} if there is a morphism $H\colon IA\r B$ with $Hi=(f,g)$. The maps $f,g$ are \emph{homotopic} $f\simeq g$ if there exists a weak equivalence $h\colon B\st{\sim}\r B'$ such that $hf$ and $hg$ are strictly homotopic. `Being homotopic' is a natural equivalence relation and the quotient category is denoted by $\pi\C{W}$. 
The homotopy category $\ho\C{W}$ is obtained by calculus of left fractions in $\pi\C{W}$. Objects in $\ho\C{W}$ are the same as in $\C{W}$.
A morphism $A\r B$ in $\ho\C{W}$ is represented by a diagram in $\C{W}$,
$$A\mathop{\To}\limits_{\alpha_1} X\mathop{\longleftarrow}\limits^\sim_{\alpha_2} B.$$
Another diagram
$$A\mathop{\To}\limits_{\alpha'_1} Y\mathop{\longleftarrow}\limits^\sim_{\alpha'_2} B$$
represents the same morphism if there is a diagram in $\C{W}$
$$\xymatrix{&X^{\estrella}\ar@{<-}[ld]_{\alpha_1}\ar@{<-}[rd]_{\sim}^{\alpha_2}&\\
A^{\estrella}&Z^{\estrella}\ar@{<-}[l]\ar@{<-}[r]^\sim\ar@{<-}[u]\ar@{<-}[d]&B^{\estrella}\\&Y^{\estrella}\ar@{<-}[lu]^{\alpha'_1}\ar@{<-}[ru]^\sim_{\alpha'_2}&}$$
whose projection to $\pi\C{W}$ is commutative. Notice that, by the 2 out of 3 axiom, the vertical arrows in this diagram are also weak equivalences.
The composite of two morphisms $A\mathop{\r}\limits_{\alpha} B\mathop{\r}\limits_{\beta} C$ in $\ho \C{W}$ represented by
$$A\mathop{\To}\limits_{\alpha_1} X\mathop{\longleftarrow}\limits^\sim_{\alpha_2} B \mathop{\To}\limits_{\beta_1} Y\mathop{\longleftarrow}\limits^\sim_{\beta_2} C$$
is defined as follows. If $\beta_1$ is a cofibration then the push-out
$$\xymatrix{B\ar@{}[rd]|-{\text{push}}\ar@{>->}[r]^-{\beta_1}\ar[d]^-\sim_-{\alpha_2}&Y\ar[d]_-\sim^-{\bar{\alpha}_2}\\
X\ar@{>->}[r]_-{\bar{\beta}_1}&X\cup_BY}$$
is defined, $\bar{\alpha}_2$ is a weak equivalence, and $\beta\alpha\colon A\r C$ is represented by
$$A\mathop{\To}\limits_{\bar{\beta}_1\alpha_1} X\cup_BY\mathop{\longleftarrow}\limits^\sim_{\bar{\alpha}_2\beta_2} C.$$
In general we can factor $\beta_1$ as cofibration followed by a weak equivalence
$$\beta_1\colon B\mathop{\rightarrowtail}\limits_{\beta'_1} Z\mathop{\To}\limits^\sim_{r} Y$$
such that there is a morphism $s\colon Y\st{\sim}\into Z$ with $rs=1_Y$. The diagram
$$\xymatrix{&Y^{\estrella}\ar@{<-}[ld]_{\beta_1}\ar@{<-}[rd]_{\sim}^{\beta_2}&\\
B^{\estrella}&Z^{\estrella}\ar@{<-<}[l]^{\beta'_1}\ar@{<-}[r]^\sim_{s\beta_2}\ar[u]^-{r}&C^{\estrella}}$$
commutes in $\C{W}$, so $\beta\colon B\r C$ is also represented by
$$B\mathop{\rightarrowtail}\limits_{\beta'_1} Z\mathop{\longleftarrow}\limits^\sim_{s\beta_2} C,$$
where the first arrow is a cofibration, and we can use this representative to define the composite $\beta\alpha\colon A\r C$.


The functor
$$\zeta\colon \C{W}\To\ho\C{W}$$
is the identity on objects and sends a morphism $f\colon A\r B$ to the morphism $\zeta(f)\colon A\r B$  represented by 
$$A\mathop{\To}\limits_{f} B\mathop{\longleftarrow}\limits^\sim_{1_B} B.$$
If $f\colon A\st{\sim}\r B$ is a weak equivalence then $\zeta(f)$ is an isomorphism and $\zeta(f)^{-1}$ is represented by
$$B\mathop{\To}\limits_{1_B} B\mathop{\longleftarrow}\limits^\sim_{f} A,$$
hence a morphism $\alpha\colon A\r B$ in $\ho\C{W}$ represented by
$$A\mathop{\To}\limits_{\alpha_1} X\mathop{\longleftarrow}\limits^\sim_{\alpha_2} B$$
coincides with $\zeta(\alpha_2)^{-1}\zeta(\alpha_1)=\alpha$.
\begin{rem}
If $\alpha$ above is an isomorphism in $\ho\C{W}$ then $\zeta(\alpha_1)=\zeta(\alpha_2)\alpha$ is also an isomorphism. In particular if $\C{W}$ has a saturated class of weak equivalences then $\alpha_1\colon A\st{\sim}\r X$ is necessarily a weak equivalence.
\end{rem}

\begin{rem}
For $\C{W}=\cb(\C{E})$ the category $\pi\C{W}=\hb(\C{E})$ is usually termed the \emph{bounded homotopy category}, while $\ho\C{W}=\db(\C{E})$ is called the \emph{bounded derived category} of $\C{E}$.
\end{rem}

\section{On Waldhausen and derived $K$-theory}\label{kt}

Recall that a \emph{cofiber sequence} in a Waldhausen category $\C{W}$
$$A\into B\onto B/A$$
is a push-out diagram
$$\xymatrix{A\ar@{>->}[r]\ar[d]\ar@{}[rd]|{\text{push}}&B\ar@{->>}[d]\\{0}\ar@{>->}[r]&B/A}$$
Therefore the quotient $B/A$ is only defined up to canonical isomorphism over $B$, although the notation $B/A$ is standard in the literature.

The $K$-theories we deal with in this paper are constructed by using the Waldhausen categories $S_n\C{W}$ that we now recall.

\begin{defn}
An object $A_{\bullet\bullet}$ in the category $S_n\C{W}$, $n\geq0$, is a commutative diagram in $\C{W}$
\begin{equation}\label{escalera}
\xymatrix{
&&&&A_{nn}\\
&&&&\vdots\ar[u]\\
&&A_{22}\ar[r]&\quad\cdots\quad\ar[r]&A_{2n}\ar[u]\\
&A_{11}\ar[r]&A_{12}\ar[r]\ar[u]&\quad\cdots\quad\ar[r]&A_{1n}\ar[u]\\
A_{00}\ar[r]&A_{01}\ar[r]\ar[u]&A_{02}\ar[r]\ar[u]&\quad\cdots\quad\ar[r]&A_{0n}\ar[u]}
\end{equation}
such that $A_{ii}=0$ and $A_{ij}\into A_{ik}\onto A_{jk}$ is a cofiber sequence for all $0\leq i\leq j\leq k\leq n$. 
Notice that these conditions imply that the whole diagram is determined, up to canonical isomorphism, by the sequence of $n-1$ composable cofibrations
\begin{equation}\label{cadena}
A_{01}\into A_{02}\into\cdots\into A_{0n}.
\end{equation}

A morphism $A_{\bullet\bullet}\r B_{\bullet\bullet}$ in $S_n\C{W}$ is a natural transformation between diagrams given by morphisms $A_{ij}\r B_{ij}$ in $\C{W}$. The category $S_n\C{W}$  is a Waldhausen category. A morphism $A_{\bullet\bullet}\st{\sim}\r B_{\bullet\bullet}$ is a weak equivalence if all morphisms $A_{ij}\st{\sim}\r B_{ij}$ are weak equivalences in $\C{W}$. A cofibration $A_{\bullet\bullet}\into B_{\bullet\bullet}$ is a morphism such that $A_{ij}\into B_{ij}$ and $B_{ij}\cup_{A_{ij}}A_{ik}\into B_{ik}$ are cofibrations, $0\leq i\leq j\leq k\leq n$. The distinguished zero object is the diagram with $0$ in all entries. 


The categories $S_n\C{W}$ assemble to a simplicial category $S.\C{W}$. The face functor $d_i\colon S_n\C{W}\r S_{n-1}\C{W}$ is defined by removing the $i^{\text{th}}$ row and the $i^{\text{th}}$  column, and the degeneracy functor $s_i\colon S_n\C{W}\r S_{n+1}\C{W}$ is defined by duplicating the $i^{\text{th}}$ row and the $i^{\text{th}}$ column, $0\leq i\leq n$. Faces and degeneracies are exact functors. For the definition of the simplicial structure it is crucial to consider the whole diagram \eqref{escalera} instead of just \eqref{cadena}.
\end{defn}

One can obtain a pointed space out of the simplicial category $S.\C{W}$ as follows. We restrict to the subcategories of weak equivalences $wS.\C{W}$, then we take levelwise the nerve in order to get a bisimplicial set $\ner wS.\C{W}$, we consider the diagonal simplicial set $\diag\ner wS.\C{W}$, and its geometric realization
$$|\diag\ner wS.\C{W}|.$$
This pointed space, actually a reduced $CW$-complex, is the $1$-stage of the \emph{Waldhausen $K$-theory} spectrum $K(\C{W})$ \cite{akts}, which is an $\Omega$-spectrum, hence the $K$-theory groups of $\C{W}$ are the homotopy groups 
\begin{eqnarray*}
K_n(\C{W})&=&\pi_{n+1}|\diag\ner wS.\C{W}|,\quad n\geq0.
\end{eqnarray*}

We now assume that $\C{W}$ has cylinders and satisifies the 2 out of 3 axiom, so that the associated right pointed derivator $\der\C{W}$ is defined, see \cite[Corollary 2.24 and the duals of Lemmas 4.2 and 4.3]{ciscd}. Then the Waldhausen categories $S_n\C{W}$ also have cylinders and satisfy the 2 out of 3 axiom. We will neither recall the notion of derivator nor the definition of the derivator $\der\C{W}$ but just the $K$-theory of $\der\C{W}$, we refer the interested reader to \cite{derivateurs, ktdt, sdckt1,cht}. For this we consider the homotopy categories $\ho S_n\C{W}$ and the subgroupoids of isomorphisms $i\ho S_n\C{W}$. These groupoids form a simplicial groupoid $i\ho S.\C{W}$ and we can consider the pointed space
$$|\diag\ner i\ho S.\C{W}|,$$
which is the $1$-stage of Garkusha's \emph{derived $K$-theory} $\Omega$-spectrum $DK(\C{W})$. 

Garkusha \cite{sdckt2} considers derived $K$-theory for $\C{W}=\cb(\C{E})$, and more generally for $\C{W}$ a nice complicial biWaldhausen category, although the definition immediately extends to Waldhausen categories with cylinders satisfying the 2 out of 3 axiom, as indicated here. Moreover, Garkusha shows that there is a natural weak equivalence $DK(\C{W})\st{\sim}\r K(\der\C{W})$ between the derived $K$-theory spectrum of a nice complicial biWaldhausen category $\C{W}$ and the $K$-theory spectrum of the associated derivator $\der\C{W}$, compare \cite[Corollary 4.3]{sdckt2}. Nevertheless \cite[Corollary 4.3]{sdckt2} only uses the fact that all morphisms in $\C{W}$ factor as a cofibration followed by a weak equivalence, compare also \cite[Lemmas 4.1 and 4.2]{sdckt2}, so we also have a natural weak equivalence $DK(\C{W})\st{\sim}\r K(\der\C{W})$ for $\C{W}$ a Waldhausen category with cylinders satisfying the 2 out of 3 axiom, and therefore
\begin{eqnarray*}
K_n(\der\C{W})&\cong&\pi_{n+1}|\diag\ner i\ho S.\C{W}|,\quad n\geq0.
\end{eqnarray*}

The functors $\zeta\colon S_n\C{W}\r\ho S_n\C{W}$ restrict to  $wS_n\C{W}\r i\ho S_n\C{W}$. These functors give rise to a map
$$|\diag\ner wS.\C{W}|\To |\diag\ner i\ho S.\C{W}|$$
which induces the comparison homomorphisms in homotopy groups,
$$\mu_n\colon K_n(\C{W})\To K_n(\der\C{W}),\quad n\geq0.$$
This map is actually the $1$-stage of a comparison map of spectra
\begin{equation}\label{cms}
K(W)\To K(\der\C{W}).
\end{equation}




In the rest of this paper we will be mainly concerned with the structure of the bisimplicial sets $X=\ner wS.\C{W}$ and $Y=\ner i\ho S.\C{W}$ in low dimensions, that we now review more thoroughly.

A \emph{bisimplicial set} $Z$ consists of sets $Z_{m,n}$, $m,n\geq 0$, together with horizontal and vertical face and degeneracy maps
$$\begin{array}{cc}
d^h_i\colon Z_{m,n}\To Z_{m-1,n},&s^h_i\colon Z_{m,n}\To Z_{m+1,n},\quad 0\leq i\leq m,\\
d^v_j\colon Z_{m,n}\To Z_{m,n-1},&s^v_j\colon Z_{m,n}\To Z_{m,n+1},\quad 0\leq j\leq n,
\end{array}$$
satisfying some relations that we will not recall here, compare \cite{gj}. An element $z_{m,n}\in Z_{m,n}$ is a \emph{bisimplex} of \emph{bidegree} $(m,n)$ and \emph{total degree} $m+n$.

A generic bisimplex $z_{m,n}$ of bidegree $(m,n)$ can be depicted as the product of two geometric simplices of dimensions $m$ and $n$ with vertices labelled by the product set
$$\set{0,\dots,m}\times\set{0,\dots,n},$$
see Figs. 1 and 2.
\begin{figure}[h]\label{12gen}
\[\entrymodifiers={+0}
\begin{array}{c}\xymatrix@!R=3pc@!C=3pc{\ar@{-}[r]^(0){(0,0)}^(1){(1,0)}&}\end{array}\; z_{1,0}
\quad
\begin{array}{c}
\xymatrix@!=2pc{\ar@{-}[r]^(0){(0,1)}^(1){(1,1)}\ar@{-}[d]&\ar@{-}[d]\\
\ar@{-}[r]_(0){(0,0)}_(1){(1,0)}&}
\end{array}\; z_{1,1}
\quad
\begin{array}{c}
\xymatrix@!R=3.5pc@!C=3.5pc{&\\
\ar@{-}[r]_(0){(0,0)}_(1){(1,0)}\ar@{-}[ru]^(1){(2,0)}
&\ar@{-}[u]}
\end{array}\; z_{2,0}
\]
\caption{Bisimplices of total degree $1$ and $2$.}
\end{figure}
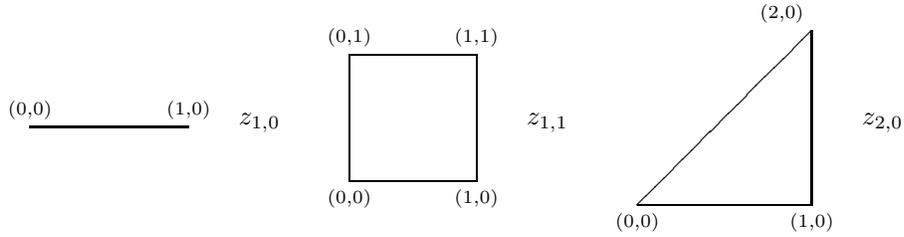
The horizontal $i^{\mathrm{th}}$ face $d^h_i(z_{m,n})$ is the face obtained by removing the interior, the vertices $(i,j)$, for all $j$, and the incident faces of the boundary. Similarly the vertical $j^{\mathrm{th}}$ face $d^v_j(z_{m,n})$ is obtained by removing the interior, the vertices $(i,j)$, for all $i$, and the incident faces of the boundary.
\begin{figure}[h]\label{3gen}
\[\entrymodifiers={+0}\xymatrix@!R=2pc@!C=0.3pc{
 \arr@{-}[rrr]^(0){(0,2)}^(1){(1,2)}&&&\\
&\arr@{.}[lu]&&&\arr@{-}[lu]\\
\arr@{-}[uu]\arr@{.}[ru]_(1){(0,1)}\arr@{-}[rrr]_(0){(0,0)}_(1){(1,0)}&&&\arr@{-}[ru]_(1){(1,1)}\arr@{-}[uu]
\arr@{.}"2,2";"2,5"|!{"1,4";"3,4"}\hole
}
\!\!\!\! z_{1,2}\qquad\qquad\quad  z_{2,1}\;\;
\xymatrix@!C=1pc@!R=0.7pc{
\arr@{-}[rrr]^(0){(0,1)}^(1){(2,1)}
\arr@{-}[rd]&&&\\&
\arr@{-}[urr]_(0){\;\;\;(1,1)}\\\arr@{-}[uu]
\arr@{.}[rrr]|(0.3333333){\hole}
\arr@{-}[rd]_(0){(0,0)}&&&\arr@{-}[uu]\\&\arr@{-}[uu]
\arr@{-}[urr]_(0){(1,0)}_(1){(2,0)}
}
\]
\vspace*{-23ex}
\[
\def\tetraco#1#2#3#4#5{{\entrymodifiers={+0}\xymatrix@!R=3.5pc@!C=0.8pc{
#5\\#5\\#5\\ z_{3,0} #5\\#5
    \arr@{-}"#1";"#2"^(0){(0,0)\qquad\qquad}
    \arr@{-}"#1";"#4"
    \arr@{-}"#2";"#3"_(0){(1,0)}
    \arr@{-}"#2";"#4"
    \arr@{-}"#3";"#4"_(0){(2,0)}_(1){(3,0)}
    \arr@{.}"#1";"#3"
    |!{"#2";"#4"}\hole
}}}
\tetraco{3,1}{4,4}{3,6}{1,4}{&&&&&}\qquad\qquad
\]
\caption{Bisimplices of total degree $3$.}
\end{figure}
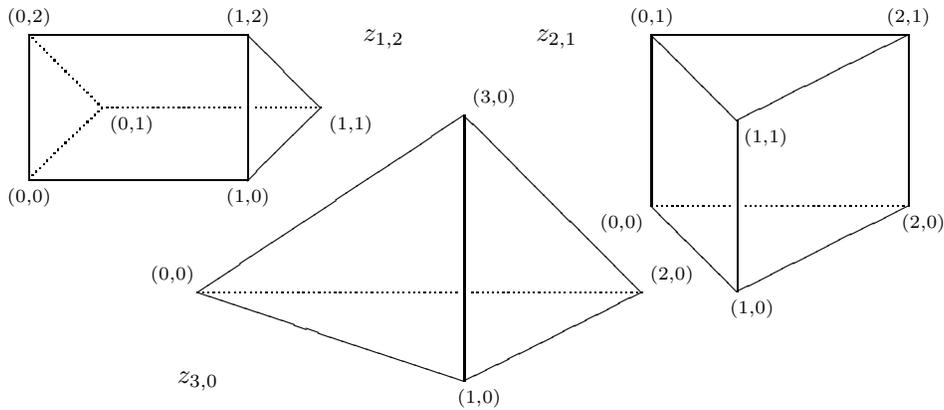

The bisimplicial sets $X$ and $Y$ are horizontally reduced, i.e. $X_{0,n}=Y_{0,n}$ are singletons for all $n\geq 0$, $X_{1,0}=Y_{1,0}$ is the set of objects in $\C{W}$, 
$X_{1,1}$ is the set of weak equivalences in $\C{W}$, $Y_{1,1}$ is the set of isomorphisms in $\ho\C{W}$, and $X_{2,0}=Y_{2,0}$ is the set of cofiber sequences, see Fig. 3.
\begin{figure}[h]
\[\entrymodifiers={+0}
\begin{array}{c}\xymatrix@!R=3pc@!C=3pc{\ar@{-}[r]|*+{\scriptstyle A^{\estrella}}&}\end{array}\quad x_{1,0}=y_{1,0}
\qquad\qquad\qquad\qquad
\begin{array}{c}
\xymatrix@!=2pc{\ar@{-}[r]|*+{\scriptstyle \bak{A}^{\estrella}}="a"\ar@{-}[d]&\ar@{-}[d]\\
\ar@{-}[r]|*+{\scriptstyle A^{\estrella}}="b"&\ar@{->}"b";"a"_{\sim}}
\end{array}\quad x_{1,1}
\]
\[\entrymodifiers={+0}
\begin{array}{c}
\xymatrix@!=2pc{\ar@{-}[r]|*+{\scriptstyle \bak{A}^{\estrella}}="a"\ar@{-}[d]&\ar@{-}[d]\\
\ar@{-}[r]|*+{\scriptstyle A^{\estrella}}="b"&\ar@{->}"b";"a"_{\cong}}
\end{array}\quad y_{1,1}
\qquad\qquad\qquad
\begin{array}{c}
\xymatrix@!R=3.5pc@!C=3.5pc{&\\
\ar@{-}[r]|*+{\scriptstyle A^{\estrella}}="a"\ar@{-}[ru]|*+{\scriptstyle B^{\estrella}}="b"
&\ar@{-}[u]|*+{\scriptstyle B^{\estrella} /A^{\estrella}}="c"\ar@{>->}"a";"b"\ar@{->>}"b";"c"}
\end{array}\quad x_{2,0}=y_{2,0}
\]
\caption{Bisimplices of total degree $1$ and $2$ in $X$ and $Y$.}
\end{figure}
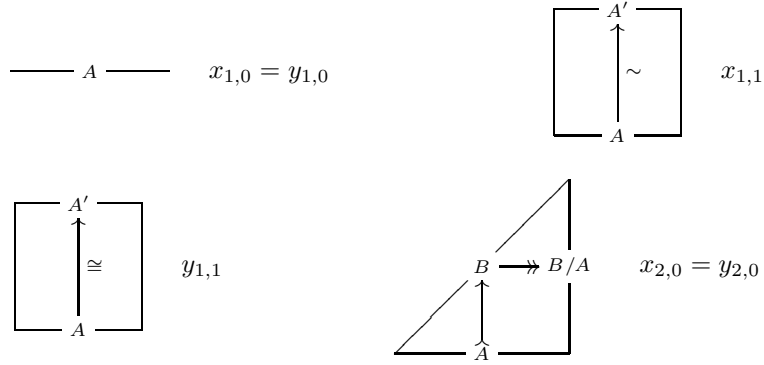 

The set $X_{1,2}$ consists of pairs of composable weak equivalences in $\C{W}$, $Y_{1,2}$ is the set of composable isomorphisms in $\ho\C{W}$, $X_{2,1}$ is the set of weak equivalences between cofiber sequences
i.e. weak equivalences in $S_2\C{W}$ which are commutative diagrams in $\C{W}$
\begin{equation}\label{wecs}
\xymatrix{\bak{A}^{\estrella}\ar@{>->}[r]& \bak{B}^{\estrella}\ar@{->>}[r]& \bak{B}^{\estrella}/\bak{A}^{\estrella}\\
A^{\estrella}\ar@{>->}[r]\ar[u]^\sim& B^{\estrella}\ar@{->>}[r]\ar[u]^\sim& B^{\estrella}/A^{\estrella}\ar[u]^\sim}
\end{equation}  
$Y_{2,1}$ is the set of isomorphisms in $\ho S_2\C{W}$,
and $X_{3,0}=Y_{3,0}$ is the set of four cofiber sequences assicated to  pairs of composable cofibrations
\begin{equation}\label{comcof}
\xymatrix{&&C^{\estrella}/B^{\estrella}\\&B^{\estrella}/A^{\estrella}\;\ar@{>->}[r]&C^{\estrella}/A^{\estrella}\ar@{->>}[u]\\A^{\estrella}\;\ar@{>->}[r]&B^{\estrella}\;\ar@{>->}[r]\ar@{->>}[u]&C^{\estrella}\ar@{->>}[u]}
\end{equation}
see Fig. 4.
\begin{figure}[h]
\[\entrymodifiers={+0}\xymatrix@!R=2pc@!C=0.3pc{
 \arr@{-}[rrr]|*+{\scriptscriptstyle C^{\estrella}}="c"&&&\\
&\arr@{.}[lu]&&&\arr@{-}[lu]\\
\arr@{-}[uu]\arr@{.}[ru]\arr@{-}[rrr]|*+{\scriptscriptstyle A^{\estrella}}="a"&&&\arr@{-}[ru]\arr@{-}[uu]
\arr@{.}"2,2";"2,5"|!{"c";"a"}\hole|*+{\scriptscriptstyle B^{\estrella}}="b"|!{"1,4";"3,4"}\hole
\arr@*{}@{.>}"a";"b"_*{\scriptscriptstyle\!\!\!\sim}
\arr@*{}@{.>}"b";"c"_*{\scriptscriptstyle\!\!\!\sim}
\arr@*{}@{->}"a";"c"^(0.3)*{\scriptscriptstyle\sim\!}
}
\!\!\!\! x_{1,2}\qquad\qquad\qquad\qquad\quad 
\xymatrix@!R=2pc@!C=0.3pc{
 \arr@{-}[rrr]|*+{\scriptscriptstyle C^{\estrella}}="c"&&&\\
&\arr@{.}[lu]&&&\arr@{-}[lu]\\
\arr@{-}[uu]\arr@{.}[ru]\arr@{-}[rrr]|*+{\scriptscriptstyle A^{\estrella}}="a"&&&\arr@{-}[ru]\arr@{-}[uu]
\arr@{.}"2,2";"2,5"|!{"c";"a"}\hole|*+{\scriptscriptstyle B^{\estrella}}="b"|!{"1,4";"3,4"}\hole
\arr@*{}@{.>}"a";"b"_*{\scriptscriptstyle\!\!\!\cong}
\arr@*{}@{.>}"b";"c"_*{\scriptscriptstyle\!\!\!\cong}
\arr@*{}@{->}"a";"c"^(0.3)*{\scriptscriptstyle\cong\!}
}
\!\!\!\! y_{1,2}
\]

\vspace{10pt}

\[\entrymodifiers={+0} x_{2,1}\qquad
\xymatrix@!C=1pc@!R=0.7pc{
\arr@{-}[rrr]|*+{\scriptscriptstyle \bak{B}^{\estrella}}="bb"
\arr@{-}[rd]|*+{\scriptscriptstyle \bak{A}^{\estrella}}="aa"&&&\\&
\arr@{-}[urr]|*+{\scriptscriptstyle \bak{B}^{\estrella}/\bak{A}^{\estrella}}="bbaa"\\\arr@{-}[uu]
\arr@{.}[rrr]|(0.1666666){\hole}|(0.3333333){\hole}
         |*+{\scriptscriptstyle B^{\estrella}}="b"|(0.6666666){\hole}
\arr@{-}[rd]|*+{\scriptscriptstyle A^{\estrella}}="a"&&&\arr@{-}[uu]\\&\arr@{-}[uu]
\arr@{-}[urr]|*+{\scriptscriptstyle B^{\estrella}/A^{\estrella}}="ba"
\arr@*{}@{>->} "aa";"bb"
\arr@*{}@{>.>}"a";"b"|{\hole}
\arr@*{}@{->>}  "bb";"bbaa"
\arr@*{}@{.>>} "b";"ba"
\arr@*{}@{->}  "a";"aa"^*{\scriptscriptstyle\sim}
\arr@*{}@{->}  "ba";"bbaa"^*{\scriptscriptstyle\sim}
\arr@*{}@{.>} "b";"bb"^(0.3)*{\scriptscriptstyle\sim}
|!{"2,2";"bbaa"}{\hole}}
\qquad\qquad\qquad
y_{2,1}\qquad
\xymatrix@!C=1pc@!R=0.7pc{
\arr@{-}[rrr]|*+{\scriptscriptstyle \bak{B}^{\estrella}}="bb"
\arr@{-}[rd]|*+{\scriptscriptstyle \bak{A}^{\estrella}}="aa"&&&\\&
\arr@{-}[urr]|*+{\scriptscriptstyle \bak{B}^{\estrella}/\bak{A}^{\estrella}}="bbaa"\\\arr@{-}[uu]
\arr@{.}[rrr]|(0.3333333){\hole}
         |*+{\scriptscriptstyle B^{\estrella}}="b"
\arr@{-}[rd]|*+{\scriptscriptstyle A^{\estrella}}="a"&&&\arr@{-}[uu]\\&\arr@{-}[uu]
\arr@{-}[urr]|*+{\scriptscriptstyle B^{\estrella}/A^{\estrella}}="ba"
\arr@*{}@{>->} "aa";"bb"
\arr@*{}@{>.>}"a";"b"|{\hole}
\arr@*{}@{->>}  "bb";"bbaa"
\arr@*{}@{.>>} "b";"ba"
\arr@*{}@{} "a";"ba"|{}="c"
\arr@*{}@{} "aa";"bbaa"|{}="cc"
\arr@*{}@{.>} "c";"cc"_{\text{ in }\ho(S_2\C{W})}_(.6)\cong|(.12){\hole}|(.2){\hole}|(.28){\hole}|(.78){\hole}|(.85){\hole}}
\]

\[\entrymodifiers={+0} 
\def\tetra#1#2#3#4#5{{\entrymodifiers={+0}\xymatrix@!R=3.5pc@!C=0.8pc{
#5\\#5\\#5\\ x_{3,0}=y_{3,0} #5\\#5
    \arr@{-}"#1";"#2"|*+{\scriptscriptstyle A^{\estrella}}="a"
    \arr@{-}"#1";"#4"|*+{\scriptscriptstyle C^{\estrella}}="c"
    \arr@{-}"#2";"#3"|*+{\scriptscriptstyle B^{\estrella}/A^{\estrella}}="ba"
    \arr@{-}"#2";"#4"|*+{\scriptscriptstyle C^{\estrella}/A^{\estrella}}="ca"
    \arr@{-}"#3";"#4"|*+{\scriptscriptstyle C^{\estrella}/B^{\estrella}}="cb"
    \arr@{.}"#1";"#3"|!{"a";"c"}\hole|*+{\scriptscriptstyle B}="b"|!{"#2";"#4"}\hole|!{"ba";"ca"}\hole
    \arr@{.>>}@*{}"b";"ba"|!{"#2";"#4"}\hole
    \arr@{.>>}@*{}"c";"cb"|!{"#2";"#4"}\hole
    \arr@{>.>}@*{}"a";"b"
    \arr@{>->} @*{}"a";"c"
    \arr@{>.>}@*{}"b";"c"
    \arr@{>->} @*{}"ba";"ca"
    \arr@{->>} @*{}"c";"ca"
    \arr@{->>} @*{}"ca";"cb"}}}
\tetra{3,1}{4,4}{3,6}{1,4}{&&&&&}
\]
\caption{Bisimplices of total degree $3$ in $X$ and $Y$.}
\end{figure}
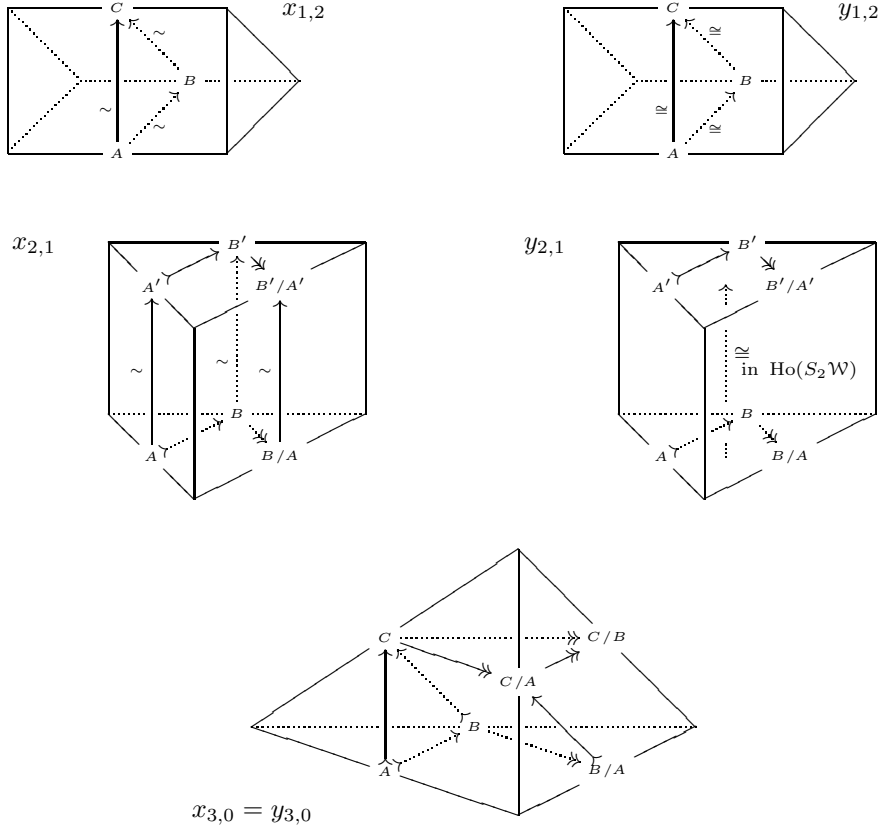 


Suppose that $\C{W}$ has a saturated class of weak equivalences. Then the categories $S_n\C{W}$ inherit this property. Therefore the isomorphism $y_{2,1}$ in $\ho S_2\C{W}$ is represented by a commutative diagram in $\C{W}$
\begin{equation}\label{wes2}
\xymatrix{\bak{A}^{\estrella}\ar@{>->}[r]& \bak{B}^{\estrella}\ar@{->>}[r]& \bak{B}^{\estrella}/\bak{A}^{\estrella}\\
X^{\estrella}\ar@{>->}[r]\ar@{<-}[u]^\sim_{\alpha_2}\ar@{<-}[d]_\sim^{\alpha_1}& Y^{\estrella}\ar@{->>}[r]\ar@{<-}[u]^\sim_{\beta_2}\ar@{<-}[d]_\sim^{\beta_1}& Y^{\estrella}/X^{\estrella}\ar@{<-}[u]^\sim_{\gamma_2}\ar@{<-}[d]_\sim^{\gamma_1}\\
A^{\estrella}\ar@{>->}[r]& B^{\estrella}\ar@{->>}[r]& B^{\estrella}/A^{\estrella}}
\end{equation}
where the horizontal lines are cofiber sequences and the vertical arrows are weak equivalences. 
The face $d_1^v(y_{2,1})$ is a cofiber sequence in $\C{W}$ which is the source of the isomorphism in $\ho S_2\C{W}$, and the face $d_0^v(y_{2,1})$ is the target. 
The faces $d_2^h(y_{2,1})$, $d_1^h(y_{2,1})$, $d_0^h(y_{2,1})$ correspond, in this order, to the isomorphisms $\alpha$, $\beta$, $\gamma$ in $\ho S_2\C{W}$ represented by the vertical lines in the previous diagram.

Notice that the representative of $y_{2,1}$ corresponds to the pasting of two bisimplices of bidegree $(2,1)$ in $X$ through a common face, see Fig. 5.
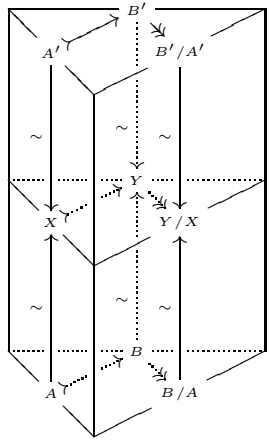
\begin{figure}[h]
\[\entrymodifiers={+0}
\xymatrix@!C=1pc@!R=0.7pc{
\arr@{-}[rrr]|*+{\scriptscriptstyle \bak{B}^{\estrella}}="bb"
\arr@{-}[rd]|*+{\scriptscriptstyle \bak{A}^{\estrella}}="aa"&&&\\&
\arr@{-}[urr]|*+{\scriptscriptstyle \bak{B}^{\estrella}/\bak{A}^{\estrella}}="bbaa"\\\arr@{-}[uu]
\arr@{.}[rrr]|(0.1666666){\hole}|(0.3333333){\hole}
         |*+{\scriptscriptstyle Y^{\estrella}}="y"|(0.6666666){\hole}
\arr@{-}[rd]|*+{\scriptscriptstyle X^{\estrella}}="x"&&&\arr@{-}[uu]\\&\arr@{-}[uu]
\arr@{-}[urr]|*+{\scriptscriptstyle Y^{\estrella}/X^{\estrella}}="yx"
\arr@*{}@{>.>}"x";"y"|{\hole}
\arr@*{}@{->>}  "bb";"bbaa"
\arr@*{}@{.>>} "y";"yx"
\arr@*{}@{<-}  "x";"aa"^*{\scriptscriptstyle\sim}
\arr@*{}@{<-}  "yx";"bbaa"^*{\scriptscriptstyle\sim}
\arr@*{}@{<.} "y";"bb"^(0.3)*{\scriptscriptstyle\sim}
|!{"2,2";"bbaa"}{\hole}\\
\arr@{-}[uu]
\arr@{.}[rrr]|(0.1666666){\hole}|(0.3333333){\hole}
         |*+{\scriptscriptstyle B^{\estrella}}="b"|(0.6666666){\hole}
\arr@{-}[rd]|*+{\scriptscriptstyle A^{\estrella}}="a"&&&\arr@{-}[uu]\\&\arr@{-}[uu]
\arr@{-}[urr]|*+{\scriptscriptstyle B^{\estrella}/A^{\estrella}}="ba"
\arr@*{}@{>->} "aa";"bb"
\arr@*{}@{>.>}"a";"b"|{\hole}
\arr@*{}@{.>>} "b";"ba"
\arr@*{}@{->}  "a";"x"^*{\scriptscriptstyle\sim}
\arr@*{}@{->}  "ba";"yx"^*{\scriptscriptstyle\sim}
\arr@*{}@{.>} "b";"y"^(0.3)*{\scriptscriptstyle\sim}|(0.63){\hole}
|!{"2,2";"bbaa"}{\hole}}
\]
\caption{A representative of $y_{2,1}$ given by the pasting of two $x_{2,1}$'s.}
\end{figure}

The degenerate bisimplices of total degree $1$ and $2$ in $X$ and $Y$ are depicted in Fig. 6.
\begin{figure}[h]
\[\entrymodifiers={+0}
\begin{array}{c}\xymatrix@!R=3pc@!C=3pc{\ar@{-}[r]|*+{\scriptstyle 0}&}\end{array}\quad s^h_0(0)
\qquad\qquad
\begin{array}{c}
\xymatrix@!=2pc{\ar@{-}[r]|*+{\scriptstyle A^{\estrella}}="a"\ar@{-}[d]&\ar@{-}[d]\\
\ar@{-}[r]|*+{\scriptstyle A^{\estrella}}="b"&\ar@{=}"b";"a"}
\end{array}\quad
s^v_0(A)
\]

\[\entrymodifiers={+0}
\begin{array}{c}
\xymatrix@!R=3pc@!C=3pc{&\\
\ar@{-}[r]|*+{\scriptstyle A^{\estrella}}="a"\ar@{-}[ru]|*+{\scriptstyle A^{\estrella}}="b"
&\ar@{-}[u]|*+{\scriptstyle 0}="c"\ar@{=}"a";"b"\ar@{->>}"b";"c"}
\end{array}\quad
s^h_1(A)
\qquad\qquad
\begin{array}{c}
\xymatrix@!R=3pc@!C=3pc{&\\
\ar@{-}[r]|*+{\scriptstyle 0}="a"\ar@{-}[ru]|*+{\scriptstyle A^{\estrella}}="b"
&\ar@{-}[u]|*+{\scriptstyle A^{\estrella}}="c"\ar@{>->}"a";"b"\ar@{=}"b";"c"}
\end{array}\quad s^h_0(A)
\]
\caption{Degenerate bisimplices of total degree $1$ and $2$ in $X$ and $Y$.}
\end{figure}
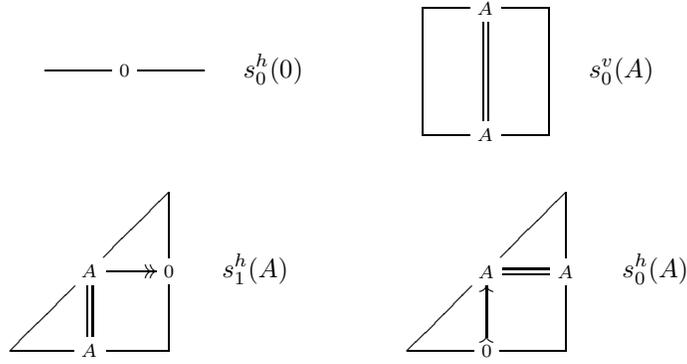

The choice of binary coproducts $A\vee B$ in $\C{W}$ gives rise to a biexact functor $\vee\colon\C{W}\times\C{W}\r\C{W}$ which induces maps of bisimplical sets \cite{akts,sdckt2}
$$\begin{array}{c}X\times X\st{\vee}\To X,\\{}\\
Y\times Y\st{\vee}\To Y,\end{array}$$
in the obvious way. 
These maps induce co-$H$-multiplications in $|\diag X|$ and $|\diag Y|$, which come from the fact that they are infinite loop spaces.

\section{Abelian $2$-groups}

In this section we recall the definition of stable quadratic modules, introduced in
\cite[Definition IV.C.1]{ch4c}. Related structures are stable crossed modules \cite{2cm} and symetric categorical groups \cite{ccscg, sccscg}. All these algebraic structures yield equivalent $2$-dimensional extensions of the theory of abelian groups. Among them stable quadratic modules are specially convenient since they are as small and strict as possible.

\begin{defn}\label{ob}
A \emph{stable quadratic module} $C_*$ is a diagram of group homomorphisms
$$C_0^{ab}\otimes C_0^{ab}\st{\grupo{\cdot,\cdot}}\longrightarrow C_1\st{\partial}\longrightarrow C_0$$
such that given $c_i,d_i\in C_i$, $i=0,1$,
\begin{enumerate}
\item $\partial\grupo{c_0,d_0}=[d_0,c_0]$,
\item $\grupo{\partial (c_1),\partial (d_1)}=[d_1,c_1]$,
\item $\grupo{c_0,d_0}+\grupo{d_0,c_0}=0$.
\end{enumerate}
Here $[x,y]=-x-y+x+y$ is the commutator of two elements
$x,y\in K$ in a group $K$, and $K^{ab}$ is the abelianization of
$K$. 

A \emph{morphism} $f\colon C_*\rightarrow D_*$ of stable quadratic modules is given by group homomorphisms $f_i\colon C_i\rightarrow D_i$, $i=0,1$, compatible with the structure homomorphisms of $C_*$ and $D_*$, i.e. $f_0\partial=\partial f_1$ and $f_1\grupo{\cdot,\cdot}=\grupo{f_0,f_0}$.
\end{defn}

\begin{rem}\label{crom}
It follows from Definition \ref{ob} that the image of $\grupo{\cdot,\cdot}$ and $\ker\partial$ are central in $C_1$, the groups $C_0$ and $C_1$ have nilpotency class $2$, and $\partial(C_1)$ is a normal subgroup of $C_0$.

There is a natural right action of $C_0$ on $C_1$ defined by
\begin{align*}
c_1^{c_0}&=c_1+\grupo{c_0,\partial(c_1)}.
\end{align*}
The axioms of a stable quadratic module imply that
commutators in $C_0$ act trivially on $C_1$, and that $C_0$ acts
trivially on the image of $\grupo{\cdot,\cdot}$ and on $\ker\partial$.

The action gives  $\partial\colon C_1\to C_0$ the structure of a crossed
module. Indeed a stable quadratic module is the same as a commutative
monoid in the category of crossed modules such that the monoid product of two
elements in $C_0$ vanishes when one of them is a commutator, see \cite[Lemma 4.18]{1tk}.
\end{rem}

\begin{rem}
The forgetful functor from stable quadratic modules to pairs of sets
$$\CC{squad}\To\CC{Set}\times\CC{Set}\colon C_*\mapsto (C_0,C_1)$$
has a left adjoint. This makes possible to define a free stable quadratic module with generating set $E_0$ in dimension $0$ and $E_1$ in dimension $1$. One can more generally define a stable quadratic module by a presentation with generators and relations in degrees $0$ and $1$. The explicit construction of a stable quadratic module with a given presentation can be found in the appendix of \cite{1tk}. For the purposes of this paper it will be enough to assume the existence of this construction, satisfying the obvious universal property as in the case of groups.
\end{rem}

We now recall the connection of stable quadratic modules with stable homotopy theory.

\begin{defn}
The \emph{homotopy groups} of a stable quadratic module $C_*$ are
\begin{eqnarray*}
\pi_0 C_*&=& C_0/\partial(C_1),\\
\pi_1 C_*&=& \ker[\partial\colon C_1\r C_0].
\end{eqnarray*}
Notice that these groups are abelian.
Homotopy groups are obviously functors in the category $\CC{squad}$ of stable quadratic modules. A morphism in $\CC{squad}$ is a \emph{weak equivalence} if it induces isomorphisms in $\pi_0$ and $\pi_1$. The \emph{$k$-invariant} of $C_*$ is the natural homomorphism
$$k\colon \pi_0C_*\otimes \mathbb{Z}/2\To \pi_1 C_*$$
defined as $k(x\otimes 1)=\grupo{x,x}$. 
\end{defn}

Weak equivalences in the Bousfield-Friedlander category $\CC{Spec}_0$ of connective spectra of simplicial sets \cite{htgamma} are also morphisms inducing isomorphisms in homotopy groups. Extending Definition \ref{hotsat}, if $\CC{C}$ is a category endowed with a class of weak equivalences we denote by $\ho\CC{C}$ the localization of $\CC{C}$ with respect to weak equivalences in the sense of \cite{gz}.

\begin{lem}{\cite[Lemma 4.22]{1tk}}\label{llave}
There is defined a functor 
$$\lambda_0\colon\ho\CC{Spec}_0\To\ho\CC{squad}$$
together with natural isomorphisms
\begin{eqnarray*}
\pi_0 \lambda_0X&\cong& \pi_0X,\\
\pi_1 \lambda_0X&\cong& \pi_1X.
\end{eqnarray*}
The $k$-invariant of $\lambda_0 X$ corresponds to the action of the stable Hopf map in the stable homotopy groups of spheres $0\neq\eta\in \pi_1^s\cong\mathbb{Z}/2$, $$\pi_0X\otimes\mathbb{Z}/2\To \pi_1X\colon x\otimes 1\mapsto x\cdot\eta.$$
Moreover, $\lambda_0$ restricts to an equivalence of categories on the full subcategory of spectra with homotopy groups concentrated in dimensions $0$ and $1$.
\end{lem}

We interpret this lemma as follows. Chain complexes of abelian groups $$\cdots\r 0\r B_1\st{\partial}\To B_0\r 0\r\cdots$$ do not model all spectra with homotopy groups concentrated in dimensions $0$ and $1$ since these complexes neglect the stable Hopf map. However these spectra are modelled by stable quadratic modules, which can be regarded as non-abelian chain complexes
$$\cdots\r 0\r C_1\st{\partial}\To C_0\r 0\r\cdots$$
endowed with an extra map
$$C_0^{ab}\otimes C_0^{ab}\st{\grupo{\cdot,\cdot}}\To C_1$$
which keeps track of the behaviour of commutators in $C_1$ and $C_0$. The homology of this non-abelian chain complex are the homotopy groups of the corresponding spectrum. Moreover, squaring the bracket $\grupo{\cdot,\cdot}$ we recover the action of the stable Hopf map.

In Section \ref{kt} we recalled that $K$-theory spectra are spectra of topological spaces. In this section we have stated Lemma \ref{llave} for spectra of simplicial sets. The geometric realization functor from simplicial sets to spaces induces an equivalence between the
the homotopy categories of spectra of simplicial sets and spectra of topological spaces. Therefore in the next section we feel free to apply the functor $\lambda_0$ in Lemma \ref{llave} to $K$-theory spectra.

\section{Algebraic models for lower $K$-theory}

In \cite{1tk} we define a stable quadratic module $\D{D}_*\C{W}$ for any Waldhausen category $\C{W}$ which is naturally isomorphic to $\lambda_0K(\C{W})$ in the homotopy category of stable quadratic modules, therefore $\D{D}_*\C{W}$ is a model for the $1$-type of the Waldhausen $K$-theory of $\C{W}$. The stable quadratic module $\D{D}_*\C{W}$ is defined by a presentation with as few generators as possible. We now recall this presentation.

\begin{defn}
We define $\D{D}_*\C{W}$ as the stable quadratic module generated in
dimension zero by the symbols
\begin{enumerate}
\renewcommand{\labelenumi}{\genn{\arabic{enumi}}}
\item $[A^{\estrella}]$ for any object in $\C{W}$, 
\end{enumerate}
and in dimension one by
\begin{enumerate}
\renewcommand{\labelenumi}{\genn{\arabic{enumi}}}
\setcounter{enumi}{1}
\item $[A^{\estrella}\st{\sim}\rightarrow \bak{A}^{\estrella}]$ for any weak equivalence, 
\item $[A^{\estrella}\into B^{\estrella}\onto B^{\estrella}/A^{\estrella}]$ for any cofiber sequence,
\end{enumerate}
such that the following relations hold.
\begin{enumerate}
\renewcommand{\labelenumi}{\reln{\arabic{enumi}}}
\item $\partial[A^{\estrella}\st{\sim}\rightarrow \bak{A}^{\estrella}]=-[\bak{A}^{\estrella}]+[A^{\estrella}]$.
\item $\partial[A^{\estrella}\into B^{\estrella}\onto B^{\estrella}/A^{\estrella}]=-[B^{\estrella}]+[B^{\estrella}/A^{\estrella}]+[A^{\estrella}]$.
\item $[0]=0$.
\item $[A^{\estrella}\st{1}\rightarrow A^{\estrella}]=0$.
\item $[A^{\estrella}\st{1}\into A^{\estrella} \onto 0]=0$, $[0\into A^{\estrella}\st{1}\onto A^{\estrella}]=0$.
\item For any pair of composable weak equivalences $ A^{\estrella}\st{\sim}\rightarrow B^{\estrella}\st{\sim}\rightarrow C^{\estrella}$,
$$[A^{\estrella}\st{\sim}\rightarrow C^{\estrella}]=[B^{\estrella}\st{\sim}\rightarrow C^{\estrella}]+[A^{\estrella}\st{\sim}\rightarrow B^{\estrella}].$$
\item For any weak equivalence between cofiber sequences in $\C{W}$, given by a commutative diagram \eqref{wecs}, 
we have
\begin{align*}
[A^{\estrella}\st{\sim}\rightarrow \bak{A}^{\estrella}]+[B^{\estrella}/A^{\estrella}\st{\cong}\rightarrow \bak{B}^{\estrella}/\bak{A}^{\estrella}]^{[A^{\estrella}]}=&
-[\bak{A}^{\estrella}\into \bak{B}^{\estrella}\onto \bak{B}^{\estrella}/\bak{A}^{\estrella}]\\&
+[B^{\estrella}\st{\sim}\rightarrow \bak{B}^{\estrella}]+[A^{\estrella}\into B^{\estrella}\onto B^{\estrella}/A^{\estrella}].
\end{align*}
\item For any commutative diagram consisting of four cofiber sequences in $\C{W}$ associated to a pair of composable cofibrations \eqref{comcof} 
we have
\begin{align*}
&[B^{\estrella}\into C^{\estrella}\onto C^{\estrella}/B^{\estrella}]+[A^{\estrella}\into
B^{\estrella}\onto B^{\estrella}/A^{\estrella}]
\\ &\quad
=[A^{\estrella}\into C^{\estrella}\onto C^{\estrella}/A^{\estrella}]+[B^{\estrella}/A^{\estrella}\into C^{\estrella}/A^{\estrella}\onto C^{\estrella}/B^{\estrella}]^{[A^{\estrella}]}.
\end{align*}
\item For any pair of objects $A^{\estrella}, B^{\estrella}$ in $\C{W}$
$$\grupo{[A],[B]}=
-[B^{\estrella}\st{i_2}\into A^{\estrella}\vee B^{\estrella}\st{p_1}\onto A^{\estrella}]
+[A^{\estrella}\st{i_1}\into A^{\estrella}\vee B^{\estrella}\st{p_2}\onto B^{\estrella}]
.$$
Here $$\xymatrix{A\ar@<.5ex>[r]^-{i_1}&A\vee
B\ar@<.5ex>[l]^-{p_1}\ar@<-.5ex>[r]_-{p_2}&B\ar@<-.5ex>[l]_-{i_2}}$$
are the inclusions and projections of a
coproduct
in $\C{W}$. 
\end{enumerate}
\end{defn}

\begin{rem}\label{semip}
The presentation of the stable quadratic module $\D{D}_*\C{W}$ is completely determined by the bisimplicial structure of $X=\ner wS.\C{W}$ and the map $\vee\colon X\times X\r X$ in total degree $\leq 3$, see Section \ref{kt}.

More precisely, $\D{D}_*\C{W}$ is generated in degree $0$ by the bisimplices of total degree $1$  and in degree $1$ by the bisimplices of total degree $2$, see Fig. 3. Relations (\relnn1) and (\relnn2) identify the image by $\partial$ of a degree $1$ generator with the summation, in an appropriate order, of the faces of the boundary of the corresponding bisimplex of total degree $2$, see again Fig. 3. Relations (\relnn3)--(\relnn5) say that degenerate bisimplices of total degree $1$ or $2$ are trivial in $\D{D}_*\C{W}$, see Fig. 6. Relations (\relnn6)--(\relnn8) tell us that the summation, in an appropriate order, of the faces of the boundary of a bisimplex of total degree $3$ is zero in $\D{D}_*\C{W}$, see Fig 4.  Finally (\relnn9) says that the bracket $\grupo{[A],[B]}$ coincides with
$$-[s_0^h(A)\vee s_1^h(B)]+[s_1^h(A)\vee s_0^h(B)],$$
i.e. it is obtained as follows. We first take the two possible degenerate bisimplices of bidegree $(2,0)$ associated to $A$ and $B$ in the following order.
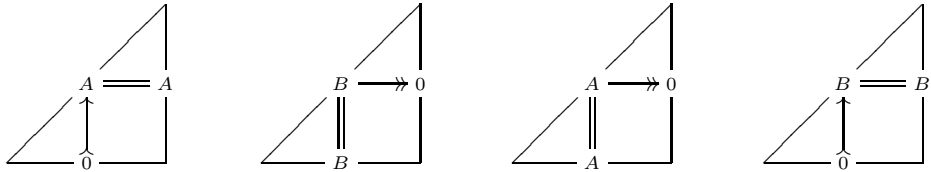
\[\entrymodifiers={+0}
\begin{array}{c}
\xymatrix@!R=3pc@!C=3pc{&\\
\ar@{-}[r]|*+{\scriptstyle 0}="a"\ar@{-}[ru]|*+{\scriptstyle A^{\estrella}}="b"
&\ar@{-}[u]|*+{\scriptstyle A^{\estrella}}="c"\ar@{>->}"a";"b"\ar@{=}"b";"c"}
\end{array}
\qquad
\begin{array}{c}
\xymatrix@!R=3pc@!C=3pc{&\\
\ar@{-}[r]|*+{\scriptstyle B^{\estrella}}="a"\ar@{-}[ru]|*+{\scriptstyle B^{\estrella}}="b"
&\ar@{-}[u]|*+{\scriptstyle 0}="c"\ar@{=}"a";"b"\ar@{->>}"b";"c"}
\end{array}
\qquad
\begin{array}{c}
\xymatrix@!R=3pc@!C=3pc{&\\
\ar@{-}[r]|*+{\scriptstyle A^{\estrella}}="a"\ar@{-}[ru]|*+{\scriptstyle A^{\estrella}}="b"
&\ar@{-}[u]|*+{\scriptstyle 0}="c"\ar@{=}"a";"b"\ar@{->>}"b";"c"}
\end{array}
\qquad
\begin{array}{c}
\xymatrix@!R=3pc@!C=3pc{&\\
\ar@{-}[r]|*+{\scriptstyle 0}="a"\ar@{-}[ru]|*+{\scriptstyle B^{\estrella}}="b"
&\ar@{-}[u]|*+{\scriptstyle B^{\estrella}}="c"\ar@{>->}"a";"b"\ar@{=}"b";"c"}
\end{array}
\]
We then take the coproduct of the first and the second pair of degenerate bisimplices.
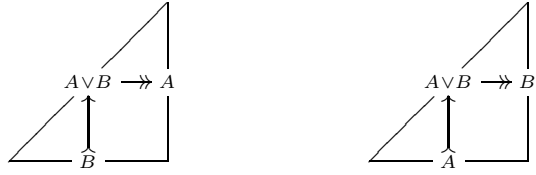
\[\entrymodifiers={+0}
\begin{array}{c}
\xymatrix@!R=3pc@!C=3pc{&\\
\ar@{-}[r]|*+{\scriptstyle B^{\estrella}}="a"\ar@{-}[ru]|*+{\scriptstyle A^{\estrella}\vee B^{\estrella}}="b"
&\ar@{-}[u]|*+{\scriptstyle A^{\estrella}}="c"\ar@{>->}"a";"b"\ar@{->>}"b";"c"}
\end{array}
\qquad\qquad\qquad
\begin{array}{c}
\xymatrix@!R=3pc@!C=3pc{&\\
\ar@{-}[r]|*+{\scriptstyle A^{\estrella}}="a"\ar@{-}[ru]|*+{\scriptstyle A^{\estrella}\vee B^{\estrella}}="b"
&\ar@{-}[u]|*+{\scriptstyle B^{\estrella}}="c"\ar@{>->}"a";"b"\ar@{->>}"b";"c"}
\end{array}
\]
Finally we take the difference between the corresponding generators in $\D{D}_1\C{W}$ 
\[\entrymodifiers={+0}
-\left[\begin{array}{c}
\xymatrix@!R=3pc@!C=3pc{&\\
\ar@{-}[r]|*+{\scriptstyle B^{\estrella}}="a"\ar@{-}[ru]|*+{\scriptstyle A^{\estrella}\vee B^{\estrella}}="b"
&\ar@{-}[u]|*+{\scriptstyle A^{\estrella}}="c"\ar@{>->}"a";"b"\ar@{->>}"b";"c"}
\end{array}\right]
\qquad+\qquad
\left[\begin{array}{c}
\xymatrix@!R=3pc@!C=3pc{&\\
\ar@{-}[r]|*+{\scriptstyle A^{\estrella}}="a"\ar@{-}[ru]|*+{\scriptstyle A^{\estrella}\vee B^{\estrella}}="b"
&\ar@{-}[u]|*+{\scriptstyle B^{\estrella}}="c"\ar@{>->}"a";"b"\ar@{->>}"b";"c"}
\end{array}\right].
\]
There is a non-abelian Eilenberg-Zilber theorem behind this formula, compare \cite[Theorem 4.10 and Example 4.13]{1tk}.
\end{rem}

The main result of \cite{1tk} is the following theorem.

\begin{thm}{\cite[Theorem 1.7]{1tk}}\label{main1tk}
Let $\C{W}$ be a Waldhausen category. There is a natural isomorphism in $\ho\CC{squad}$
$$\D{D}_*\C{W}\st{\cong}\To\lambda_0K(\C{W}).$$
\end{thm}

This result is meaningful since $\lambda_0K(\C{W})$ is huge compared with $\D{D}_*\C{W}$, while $\D{D}_*\C{W}$ is directly defined in terms of the basic structure of the Waldhausen category $\C{W}$. As a consequence we have an exact sequence of groups
$$K_1(\C{W})\hookrightarrow\D{D}_1\C{W}\st{\partial}\To \D{D}_0\C{W}\twoheadrightarrow K_0(\C{W}).$$

We now extend Theorem \ref{main1tk} to derived $K$-theory.

\begin{defn}
We define $\DD{D}_*\C{W}$ as the stable quadratic module generated in
dimension zero by the symbols
\begin{enumerate}
\renewcommand{\labelenumi}{\gens{\arabic{enumi}}}
\item $[A]$ for any object in $\C{W}$, i.e. the same as (\gennn1),
\end{enumerate}
and in dimension one by
\begin{enumerate}
\renewcommand{\labelenumi}{\gens{\arabic{enumi}}}
\setcounter{enumi}{1}
\item $[A\st{\cong}\rightarrow A']$ for any isomorphism in $\ho \C{W}$,
\item $[A\into B\onto B/A]$ for any cofiber
sequence in $\C{W}$, i.e. the same as (\gennn3),
\end{enumerate}
such that the following relations hold.
\begin{enumerate}
\renewcommand{\labelenumi}{\rels{\arabic{enumi}}}
\item $\partial[A\st{\cong}\rightarrow A']=-[A']+[A]$.
\item $=$ (\relnn2). 
\item $=$ (\relnn3). 
\item $[A\st{1}\rightarrow A]=0$.
\item $=$ (\relnn5).
\item For any pair of composable isomorphisms $A\st{\cong}\rightarrow B\st{\cong}\rightarrow C$ in $\ho\C{W}$,
$$[A\st{\cong}\rightarrow C]=[B\st{\cong}\rightarrow C]+[A\st{\cong}\rightarrow B].$$
\item For any commutative diagram in $\C{W}$ as \eqref{wes2} 
we have
\begin{align*}
[\alpha\colon A\st{\cong}\rightarrow A']+[\gamma\colon B/A\st{\cong}\rightarrow B'/A']^{[A]}=&
-[A'\into B'\onto B'/A']\\&
+[\beta\colon B\st{\cong}\rightarrow B']   
+[A\into B\onto B/A].
\end{align*}
Here $\alpha=\zeta(\alpha_2)^{-1}\zeta(\alpha_1)$, $\beta=\zeta(\beta_2)^{-1}\zeta(\beta_1)$ and $\gamma=\zeta(\gamma_2)^{-1}\zeta(\gamma_1)$.
\item $=$ (\relnn8).  
\item $=$ (\relnn9). 
\end{enumerate}
\end{defn}

If $\C{W}$ is a Waldhausen category with cylinders and a saturated class of weak equivalences then the presentation of the stable quadratic module $\DD{D}_*\C{W}$ is determined by the bisimplicial structure of $Y=\ner i\ho S.\C{W}$ and the map $\vee\colon Y\times Y\r Y$ in total degree $\leq 3$, see Section \ref{kt}, exactly in the same way as  $\D{D}_*\C{W}$ is determined by $X=\ner wS.\C{W}$ and $\vee\colon X\times X\r X$, see Remark \ref{semip}. Therefore replacing $X$ by $Y$ in the proof of \cite[Theorem 1.7]{1tk} we obtain the following result.

\begin{thm}\label{main3}
Let $\C{W}$ be a Waldhausen category with cylinders and a saturated class of weak equivalences. There is a natural isomorphism in $\ho\CC{squad}$
$$\DD{D}_*\C{W}\st{\cong}\To\lambda_0K(\der\C{W}).$$
\end{thm}

As a consequence we have an exact sequence of groups
$$K_1(\der{\C{W}})\hookrightarrow\DD{D}_1\C{W}\st{\partial}\To \DD{D}_0\C{W}\twoheadrightarrow K_0(\der{\C{W}}).$$

\begin{rem}\label{clave}
As one can easily check, taking $\lambda_0$ in the comparison map of spectra \eqref{cms} which induces $\mu_n\colon K_n(\C{W})\r K_n(\der{\C{W}})$ in homotopy groups corresponds to the natural morphism in $\CC{squad}$,
$$\bar{\mu}\colon\D{D}_*\C{W}\To\DD{D}_*\C{W},\;\;$$
$$\begin{array}{rcl}
{[}A{]}&\mapsto& {[}A{]},\\
{[}f\colon A\st{\sim}\r A'{]}&\mapsto& {[}\zeta(f)\colon A\st{\cong}\r A'{]},\\
{[}A\into B\onto B/A{]}&\mapsto& {[}A\into B\onto B/A{]},
\end{array}$$
under the natural isomorphisms of Theorems \ref{main1tk} and \ref{main3}. In particular taking $\pi_0$ and $\pi_1$ in this morphism of stable quadratic modules we obtain $\mu_0$ and $\mu_1$, respectively. This fact will be used below in the proof of Theorem \ref{main2}.
\end{rem}

\section{Proof of Theorem \ref{main2}}

Theorem \ref{main2} is a consequence of the following result.

\begin{thm}\label{el}
Let $\C{W}$ be a Waldhausen category with cylinders and a saturated class of weak equivalences. The natural morphism in $\CC{squad}$
$$\bar{\mu}\colon\D{D}_*\C{W}\To\DD{D}_*\C{W},$$
defined in Remark \ref{clave}, is an isomorphism.
\end{thm}

The key for the proof of Theorem \ref{el} is the following lemma.

\begin{lem}\label{la}
Let $\C{W}$ be a Waldhausen category with cylinders satisfying the 2 out of 3 axiom. Two weak equivalences $f,g\colon A\st{\sim}\r A'$ which are homotopic $f\simeq g$ represent the same element in $\D{D}_1\C{W}$,
$$[f\colon A\st{\sim}\r A']\;=\;[g\colon A\st{\sim}\r A'].$$
\end{lem}

\begin{proof}
Let $IA$ be a cylinder of $A$ and
$$A\vee A\mathop{\rightarrowtail}\limits_i IA\mathop{\To}\limits^\sim_p A$$
a factorization of the folding map, i.e. if $i=(i_0,i_1)$ then $pi_0=pi_1=1_A$. Since both $p$ and $1_A$ are weak equivalences we deduce from the 2 out of 3 axiom that $i_0$ and $i_1$ are also weak equivalences. Moreover, for $j=0,1$,
\begin{eqnarray*}
0&\st{\text{\scriptsize (\relnn4)}}=&[A\st{1_A}\r A]\\
&=&[pi_j\colon A\st{\sim}\r A]\\
\text{\scriptsize (\relnn6)}\quad&=&[p\colon IA\st{\sim}\r A]+[i_j\colon A\st{\sim}\r IA],
\end{eqnarray*}
therefore
$$[i_0\colon A\st{\sim}\r IA]\;\;=\;\; -[p\colon IA\st{\sim}\r A]  \;\;=\;\;[i_1\colon A\st{\sim}\r IA].$$

Furthermore, $f\simeq g$, so there is a weak equivalence $h\colon A'\st{\sim}\r A''$ and a morphism $H\colon IA\r A''$ such that $Hi_0=hf$ and $Hi_1=hg$. Again by the 2 out of 3 axiom $H$ is a weak equivalence, and
\begin{eqnarray*}
{[}h\colon A'\st{\sim}\r A''{]}+{[}f\colon A\st{\sim}\r A'{]}&\st{\text{\scriptsize (\relnn6)}}=&{[}hf=Hi_0\colon A\st{\sim}\r A''{]}\\
\text{\scriptsize(\relnn6)}\quad&=&{[}H\colon IA\st{\sim}\r A''{]}+{[}i_0\colon A\st{\sim}\r IA{]}\\
&=&{[}H\colon IA\st{\sim}\r A''{]}+{[}i_1\colon A\st{\sim}\r IA{]}\\\text{\scriptsize (\relnn6)}\quad&=&
{[}hg=Hi_1\colon A\st{\sim}\r A''{]}\\
\text{\scriptsize(\relnn6)}\quad&=&
{[}h\colon A'\st{\sim}\r A''{]}+{[}g\colon A\st{\sim}\r A'{]},
\end{eqnarray*}
hence we are done.
\end{proof}

We are now ready to prove Theorem \ref{el}.

\begin{proof}[Proof of Theorem \ref{el}]
We are going to define the inverse of $\bar{\mu}$,
$$\bar{\nu}\colon\DD{D}_*\C{W}\To\D{D}_*\C{W}.$$

We first show that
\begin{eqnarray*}
\bar{\nu}_0[A]&=&[A],\\
\bar{\nu}_1[\alpha\colon A\st{\cong}\r A']&=&-[\alpha_2\colon A'\st{\sim}\r X]+[\alpha_1\colon A\st{\sim}\r X],\\
\bar{\nu}_1{[}A\into B\onto B/A{]}&=&{[}A\into B\onto B/A{]},
\end{eqnarray*}
defines a stable quadratic module morphism $\bar{\nu}$. Here
$$A\mathop{\To}\limits^\sim_{\alpha_1} X\mathop{\longleftarrow}\limits^\sim_{\alpha_2} A'$$
is a representative of the isomorphism $\alpha$.
For this we are going to prove that the image of $[\alpha\colon A\st{\cong}\r A']$ does not depend on the choice of a representative. 

Suppose that
$$A\mathop{\To}\limits_{\alpha'_1}^\sim Y\mathop{\longleftarrow}\limits^\sim_{\alpha'_2} A'$$
also represents $\alpha$. Then there is a diagram in $\C{W}$
$$\xymatrix{&X^{\estrella}\ar@{<-}[ld]_{\alpha_1}\ar@{<-}[rd]^{\alpha_2}&\\
A^{\estrella}&Z^{\estrella}\ar@{<-}[l]|{f_1}\ar@{<-}[r]|{f_2}\ar@{<-}[u]|{g}\ar@{<-}[d]|{g'}&A'^{\estrella}\\&Y^{\estrella}\ar@{<-}[lu]^{\alpha'_1}\ar@{<-}[ru]_{\alpha'_2}&}$$
where all arrows are weak equivalences and the four triangles commute up to homotopy, so
\begin{eqnarray*}
-[\alpha_2\colon A'\st{\sim}\r X]+[\alpha_1\colon A\st{\sim}\r X]&=&-[\alpha_2\colon A'\st{\sim}\r X]-[g\colon X\st{\sim}\r Z]\\&&+[g\colon X\st{\sim}\r Z]+[\alpha_1\colon A\st{\sim}\r X]\\
\text{\scriptsize (\relnn6)}\quad&=&-[g\alpha_2\colon A'\st{\sim}\r Z]+[g\alpha_1\colon A\st{\sim}\r Z]\\
\text{\scriptsize Lemma \ref{la}}\quad&=&-[f_2\colon A'\st{\sim}\r Z]+[f_1\colon A\st{\sim}\r Z]\\
\text{\scriptsize Lemma \ref{la}}\quad&=&-[g'\alpha'_2\colon A'\st{\sim}\r Z]+[g'\alpha'_1\colon A\st{\sim}\r Z]\\
\text{\scriptsize (\relnn6)}\quad&=&-[\alpha'_2\colon A'\st{\sim}\r Y]-[g'\colon Y\st{\sim}\r Z]\\&&+[g'\colon Y\st{\sim}\r Z]+[\alpha_1'\colon A\st{\sim}\r Y]\\
&=&-[\alpha_2'\colon A'\st{\sim}\r X]+[\alpha_1'\colon A\st{\sim}\r X].
\end{eqnarray*}

Now we check that the definition of $\bar{\nu}$ on generators is compatible with the defining relations. The only non-trivial part concerns relations (\relss1), (\relss6) and  (\relss7).
Compatibility with (\relss1) follows from
\begin{eqnarray*}
\bar{\nu}_0\partial[\alpha\colon A\st{\cong}\r A']&=&\partial\bar{\nu}_1[\alpha\colon A\st{\cong}\r A']\\&=&-\partial[\alpha_2\colon A'\st{\sim}\r X]+\partial[\alpha_1\colon A\st{\sim}\r X]\\
{\scriptstyle (\relnn1)}\quad &=&-(-[X]+[A'])+(-[X]+[A])\\
&=&-[A']+[A]\\
&=&-\bar{\nu}_0[A']+\bar{\nu}_0[A].
\end{eqnarray*}
In order to check compatibility with (\relss6) we consider two composable isomorphisms in $\ho\C{W}$ $$A\mathop{\To}\limits^{\cong}_{\alpha} B\mathop{\To}\limits^{\cong}_{\beta} C$$ and we take representatives of $\alpha$, $\beta$ and $\beta\alpha$ as in the following commutative diagram of weak equivalences in $\C{W}$
$$\xymatrix@C=15pt{&&X\cup_BY\ar@{<-}[rd]^{\bar{\alpha}_2}\ar@{<-<}[ld]_{\bar{\beta}_1}&&\\&X\ar@{<-}[rd]_{\alpha_2}\ar@{<-}[ld]_{\alpha_1}\ar@{}[rr]|{\text{push}}&&Y\ar@{<-}[rd]^{\beta_2}&\\A&&B\ar@{>->}[ru]_{\beta_1}&&C}$$
Then 
\begin{eqnarray*}
\bar{\nu}_1[\beta\alpha\colon A\st{\cong}\r C]
&=&-[\bar{\alpha}_2\beta_2\colon C\st{\sim}\r X\cup_BY]+[\bar{\beta}_1\alpha_1\colon A\st{\sim}\r X\cup_BY]\\
&=&-[\bar{\alpha}_2\beta_2\colon C\st{\sim}\r X\cup_BY]+[\bar{\alpha}_2\beta_1\colon B\st{\sim}\r X\cup_BY]\\
&&-[\bar{\beta}_1\alpha_2\colon B\st{\sim}\r X\cup_BY]+[\bar{\beta}_1\alpha_1\colon A\st{\sim}\r X\cup_BY]\\
{\scriptstyle (\relnn6)}\qquad&=&-([\bar{\alpha}_2\colon Y\st{\sim}\r X\cup_BY]+[\beta_2\colon C\st{\sim}\r Y])\\
&&+[\bar{\alpha}_2\colon Y\st{\sim}\r X\cup_BY]+[\beta_1\colon B\st{\sim}\r Y]\\
&&-([\bar{\beta}_1\colon X\st{\sim}\r X\cup_BY]+[\alpha_2\colon B\st{\sim}\r X])\\
&&+[\bar{\beta}_1\colon X\st{\sim}\r X\cup_BY]+[\alpha_1\colon A\st{\sim}\r X]\\
&=&-[\beta_2\colon C\st{\sim}\r Y]+[\beta_1\colon B\st{\sim}\r Y]\\
&&-[\alpha_2\colon B\st{\sim}\r X]+[\alpha_1\colon A\st{\sim}\r X]\\
&=&\bar{\nu}_1[\beta\colon B\st{\cong}\r C]+\bar{\nu}_1[\alpha\colon A\st{\cong}\r B].
\end{eqnarray*}
Let us now check compatibility with (\relss7).
\begin{eqnarray*}
-\bar{\nu}_1[A'\into B'\onto B'/A']&&\\
+\bar{\nu}_1[\beta\colon B\st{\cong}\r B']&&\\
+\bar{\nu}_1[A\into B\onto B/A]&=&-[A'\into B'\onto B'/A']-[\beta_2\colon B'\st{\sim}\r Y]\\
&&+[X\into Y\onto Y/X]-[X\into Y\onto Y/X]\\
&&+[\beta_1\colon B\st{\sim}\r Y]+[A\into B\onto B/A]\\
{\scriptstyle (\relnn7)}\quad&=&-([\alpha_2\colon A'\st{\sim}\r X]+[\gamma_2\colon B'/A'\st{\sim}\r Y/X]^{[A']})\\
&&+[\alpha_1\colon A\st{\sim}\r X]+[\gamma_1\colon B/A\st{\sim}\r Y/X]^{[A]}\\
\text{\scriptsize Rem. \ref{crom}}\quad&=&-[\gamma_2\colon B'/A'\st{\sim}\r Y/X]-[\alpha_2\colon A'\st{\sim}\r X]\\
&&+[\alpha_1\colon A\st{\sim}\r X]+[\gamma_1\colon B/A\st{\sim}\r Y/X]\\
&&-\grupo{[A'],\partial[\gamma_2]}+\grupo{[A],\partial[\gamma_1]}\\
\text{\scriptsize Defn. \ref{ob} (2) and Rem. \ref{crom}}\quad&=&-[\alpha_2\colon A'\st{\sim}\r X]+[\alpha_1\colon A\st{\sim}\r X]\\
&&-[\gamma_2\colon B'/A'\st{\sim}\r Y/X]+[\gamma_1\colon B/A\st{\sim}\r Y/X]\\
&&+\grupo{-\partial[\alpha_2]+\partial[\alpha_1],-\partial[\gamma_2]}   
-\grupo{[A'],\partial[\gamma_2]}+\grupo{[A],\partial[\gamma_1]}\\
\text{\scriptsize (\relnn1)}\quad&=&-[\alpha_2\colon A'\st{\sim}\r X]+[\alpha_1\colon A\st{\sim}\r X]\\
&&-[\gamma_2\colon B'/A'\st{\sim}\r Y/X]+[\gamma_1\colon B/A\st{\sim}\r Y/X]\\
&&+\grupo{-(-[X]+[A'])+(-[X]+[A]),-\partial[\gamma_2]}   \\&&
+\grupo{[A'],-\partial[\gamma_2]}+\grupo{[A],\partial[\gamma_1]}\\
&=&-[\alpha_2\colon A'\st{\sim}\r X]+[\alpha_1\colon A\st{\sim}\r X]\\
&&-[\gamma_2\colon B'/A'\st{\sim}\r Y/X]+[\gamma_1\colon B/A\st{\sim}\r Y/X]\\
&&+\grupo{[A],-\partial[\gamma_2]}+\grupo{[A],\partial[\gamma_1]}\\
&=&+\bar{\nu}_1[\alpha\colon A\st{\cong}\r A']+\bar{\nu}_1[\gamma\colon B/A\st{\cong}\r B'/A']\\
&&+\grupo{\bar{\nu}_0[A],\partial\bar{\nu}_1[\gamma\colon B/A\st{\cong}\r B'/A']}\\
\text{\scriptsize Rem. \ref{crom}}\quad&=&\bar{\nu}_1[\alpha\colon A\st{\cong}\r A']
+\bar{\nu}_1[\gamma\colon B/A\st{\cong}\r B'/A']^{\bar{\nu}_0[A]}.
\end{eqnarray*}
This establishes that $\bar{\nu}$ is a well defined morphism of stable quadratic modules.

Let us now check that $\bar{\mu}\bar{\nu}=1_{\DD{D}_*\C{W}}$ and $\bar{\nu}\bar{\mu}=1_{\D{D}_*\C{W}}$. Both equations are obvious on generators $(\gennn1)=(\genss1)$ and $(\gennn3)=(\genss3)$. For $(\gennn2)$
\begin{eqnarray*}
\bar{\nu}_1\bar{\mu}_1[f\colon A\st{\sim}\r A']&=&\bar{\nu}_1[\zeta(f)\colon A\st{\cong}\r A']\\
&=&-[1_{A'}\colon A'\st{\sim}\r A']+[f\colon A\st{\sim}\r A']\\
{\scriptstyle (\relnn4)}\quad&=&[f\colon A\st{\sim}\r A'].
\end{eqnarray*}
If $\alpha\colon A\st{\cong}\r A'$ is an isomorphism in $\ho\C{W}$ we have the following equation in $\DD{D}_1\C{W}$,
\begin{eqnarray*}
0&\st{\text{\scriptsize (\relss4)}}=&[A\st{1_A}\r A]\\
&=&[\alpha^{-1}\alpha\colon A\st{\cong}\r A]\\
\text{\scriptsize (\relss6)}\quad&=&[\alpha^{-1}\colon A'\st{\cong}\r A] + [\alpha\colon A\st{\cong}\r A'],
\end{eqnarray*}
so $[\alpha^{-1}\colon A'\st{\cong}\r A]=-[\alpha\colon A\st{\cong}\r A']$. Now for $(\genss2)$
\begin{eqnarray*}
\bar{\mu}_1\bar{\nu}_1[\alpha\colon A\st{\cong}\r A']&=&-\bar{\mu}_1[\alpha_2\colon A'\st{\sim}\r X]+\bar{\mu}_1[\alpha_1\colon A\st{\sim}\r X]\\
&=&-[\zeta(\alpha_2)\colon A'\st{\cong}\r X]+[\zeta(\alpha_1)\colon A\st{\cong}\r X]\\
&=&[\zeta(\alpha_2)^{-1}\colon X\st{\cong}\r A']+[\zeta(\alpha_1)\colon A\st{\cong}\r X]\\
\text{\scriptsize (\relss6)}\quad&=&[\alpha=\zeta(\alpha_2)^{-1}\zeta(\alpha_1)\colon A\st{\cong}\r A'].
\end{eqnarray*}
The proof of Theorem \ref{el} is now finished.
\end{proof}

\begin{rem}\label{last}
Let $\C{W}$ be a Waldhausen category with cylinders satisfying the 2 out of 3 axiom. We do not assume that $\C{W}$ has a saturated class of weak equivalences. However we can endow the underlying category with a new Waldhausen category structure which does have a saturated class of weak equivalences. 

We consider the Waldhausen category $\overline{\C{W}}$ with the same underlying category as $\C{W}$. Cofibrations in $\overline{\C{W}}$ are also de same as in $\C{W}$. Weak equivalences in $\overline{\C{W}}$ are the morphisms in $\C{W}$ which are mapped to isomorphisms in $\ho\C{W}$ by the canonical functor $\zeta\colon\C{W}\r\ho\C{W}$.  Therefore weak equivalences in $\C{W}$ are also weak equivalences in $\overline{\C{W}}$ but the converse need not hold. This indeed defines a Waldhausen category $\overline{\C{W}}$ with cylinders and a saturated class of weak equivalences, and the obvious exact functor $\C{W}\r\overline{\C{W}}$ induces an isomorphism on the associated derivators $\der\C{W}\cong\der\overline{\C{W}}$, 
compare \cite[dual of Proposition 3.16]{ciscd} and \cite[Theorem 6.2.2]{cht}. Hence we have a commutative diagram for $n=0,1$,
$$\xymatrix{K_n(\C{W})\ar[r]^{\mu_n}\ar[d]&K_n(\der\C{W})\ar[d]^\cong\\
K_n(\overline{\C{W}})\ar[r]_{\mu_n}^\cong&K_n(\der\overline{\C{W}})}$$
Here the lower arrow is an isomorphism by Theorem \ref{main2}. Now we can use the `fibration theorem', \cite[1.6.7]{akts} and \cite[Theorem 11]{nktdc}, to embed the morphisms $\mu_n\colon K_n(\C{W})\r K(\der\C{W})$, $n=0,1$, in an exact sequence. More precisely, let $\C{W}_0$ be the full subcategory of $\C{W}$ spanned by the objects which are isomorphic to $0$ in $\ho\C{W}$. The category $\C{W}_0$ is a Waldhausen category where a morphism is a cofibration, resp. a weak equivalence, if and only if it is a cofibration, resp. a weak equivalence, in $\C{W}$. There is an exact sequence
$$K_1(\C{W}_0)\To K_1(\C{W})\st{\mu_1}\To K_1(\der\C{W})\st{\delta}\To K_0(\C{W}_0)\To K_0(\C{W})\st{\mu_0}\To K_0(\der\C{W})\r 0.$$

The group $K_0(\C{W}_0)$ has also been considered by Weiss in \cite{hlwc}. Weiss defines the Whitehead group of $\C{W}$ as $\operatorname{Wh}(\C{W})=K_0(\C{W}_0)$. Moreover, for any morphism $f\colon A\r A'$ which becomes an isomorphism in $\ho\C{W}$ he defines the Whitehead torsion $\tau(f)\in\operatorname{Wh}(\C{W})$, which is the obstruction for $f$ to be a weak equivalence in $\C{W}$. If $f\colon A\r A$ is an endomorphism which maps to an automorphism in $\ho\C{W}$ then one can check that $$\delta[\zeta(f)\colon A\st{\cong}\r A]\quad=\quad-\tau(f),$$ therefore an automorphism in $\ho\C{W}$ comes from a weak equivalence in $\C{W}$ if and only if its class in derivator $K_1$ comes from Waldhausen $K_1$. 
\end{rem}

\bibliographystyle{amsalpha}
\providecommand{\bysame}{\leavevmode\hbox to3em{\hrulefill}\thinspace}
\providecommand{\MR}{\relax\ifhmode\unskip\space\fi MR }
\providecommand{\MRhref}[2]{%
  \href{http://www.ams.org/mathscinet-getitem?mr=#1}{#2}
}
\providecommand{\href}[2]{#2}

\end{document}